\documentclass[12pt]{article}

\usepackage{latexsym}
\usepackage{amssymb}
\usepackage[margin=2.5cm]{geometry}

\usepackage[english]{babel}
\usepackage{amsmath}
\usepackage{amsthm}
\usepackage{lineno}

\newtheorem{theorem}{Theorem}[section]
\newtheorem{lemma}[theorem]{Lemma}
\newtheorem{corollary}[theorem]{Corollary}
\newtheorem{proposition}[theorem]{Proposition}

\theoremstyle{definition}
\newtheorem{remark}[theorem]{Remark}
\newtheorem{problem}[theorem]{Problem}

\begin{document}
\date{}
\title{A note on one-sided interval edge colorings of bipartite graphs}

\author{
Carl Johan Casselgren\footnote{Department of Mathematics, 
Link\"oping University, 
SE-581 83 Link\"oping, Sweden.
{\it E-mail address:} carl.johan.casselgren@liu.se. Supported by a grant from the Swedish
Research Council (2017-05077).}
}
\date{\today}
\maketitle

\bigskip
\noindent
{\bf Abstract.}
For a bipartite graph $G$ with parts $X$ and $Y$,
an $X$-interval coloring
is a proper edge coloring of $G$
by integers such that the colors on the edges
incident to any vertex in $X$ form an interval.
Denote by $\chi'_{int}(G,X)$ the minimum $k$ such that $G$ has an
$X$-interval coloring with $k$ colors.
The author and Toft conjectured [Discrete Mathematics 339 (2016), 2628--2639]
that there is a polynomial
$P(x)$ such that if $G$ has maximum degree at most $\Delta$, then
$\chi'_{int}(G,X) \leq P(\Delta)$.
In this short note, we prove this conjecture; in fact,
we prove that a cubic polynomial suffices.
We also deduce some improved upper bounds on $\chi'_{int}(G,X)$
for bipartite graphs with small maximum degree.

\bigskip

\noindent
\small{\emph{Keywords: One-sided interval edge coloring, Interval edge coloring, bipartite
graph, edge coloring}}

\section{Introduction}

An {\em interval coloring} of a graph  is a proper edge coloring  by integers such that the colors on the edges incident to any vertex  form an interval of integers; this notion was introduced by Asratian and Kamalian  
\cite{AsratianKamalian} (available in English as 
\cite{AsratianKamalian2}), motivated by the problem of finding compact school timetables, that is, timetables such that the lectures of each teacher and each class are scheduled at consecutive periods. Hansen \cite{Hansen} suggested another scenario (first described by Jesper Bang–Jensen): a school wishes to schedule parent–teacher conferences in time slots so that every person’s conferences occur in consecutive slots. A solution exists if and only if the bipartite graph with vertices for the people and edges for the required meetings has an interval coloring.

Not every graph has an interval coloring, since a graph
$G$ with an interval coloring must have
a proper $\Delta(G)$-edge coloring \cite{AsratianKamalian2}, where
$\Delta(G)$ denotes the maximum degree of $G$.
Sevastjanov  \cite{Sevastjanov} proved that  determining whether
a bipartite graph has an interval coloring is $\mathcal{NP}$-complete.
Nevertheless, trees \cite{Hansen, AsratianKamalian},
regular and complete bipartite graphs \cite{Hansen, AsratianKamalian},
grids \cite{timetabling4}, subcubic connected Class 1 graphs 
\cite{AsratianCasselgrenPetrosyan2},
and simple outerplanar bipartite graphs \cite{timetabling3, Axenovich}
all have interval colorings.

A well-known conjecture suggests that all {\em $(a,b)$-biregular}
graphs have interval colorings
(see e.g. \cite{Hansen, JensenToft, StiebitzToft}),
where a bipartite graph is
{\em $(a,b)$-biregular}
if all vertices in one part have degree $a$ and
all vertices in the other part
have degree $b$.
By results of \cite{Hansen} and \cite{HansonLotenToft},
all $(2,b)$-biregular graphs admit interval colorings
(the latter result was also obtained independently by
Kostochka \cite{unpublished}).
In \cite{CasselgrenToft} it is proved that every $(3,6)$-biregular
graph has an interval $7$-coloring and in
\cite{CasselgrenToft2} it was proved that large families of
$(3,5)$-biregular graphs admit interval colorings.
Several sufficient conditions for a $(3,4)$-biregular 
graph to
admit an interval $6$-coloring have been obtained
\cite{AsratianCasselgren2, Pyatkin, YangLi};
however, it remains an open question whether all such
graphs have interval colorings.

	For a bipartite graph with parts $X$ and $Y$, an {\em $X$-interval coloring}
	(or {\em one-sided interval coloring})
	is a proper edge coloring
	such that the colors on the edges incident to any vertex of $X$
	form an interval of integers. This kind of edge coloring
	seems to have been first considered in 
	\cite{AsratianKamalian, AsratianKamalian2}.
	Note that a one-sided interval coloring of a bipartite graph
	has a natural interpretation as a timetable where lectures
	are scheduled at consecutive time slots either for the teachers, or for
	the classes.
	For the graph $G$ in the scenario by Hansen discussed above,
	a one-sided interval coloring of $G$ corresponds to a schedule
	where the meetings are consecutive either for the parents, or for the teachers.

	Trivially, every bipartite graph $G$ with parts $X$ and $Y$
	has an $X$-interval coloring with $|E(G)|$ colors.
	We denote by
	$\chi'_{int}(G,X)$ the smallest integer $t$ such
	that there is an $X$-interval $t$-coloring of $G$.
	Note that, in general, the problem of computing
	$\chi'_{int}(G,X)$ is  $\mathcal{NP}$-hard;
	this follows from the fact that
	determining whether a given $(3,6)$-biregular graph has an interval
	$6$-coloring is $\mathcal{NP}$-complete
	\cite{AsratianCasselgren1}.

	Kamalian \cite{Kamalian} obtained some upper bounds on $\chi'_{int}(G,X)$
	for biregular graphs.
	As noted in \cite{CasselgrenToft}, it follows from K\"onig's edge coloring
	theorem that $\chi'_{int}(G,X) \leq 2^{\Delta(G)} \Delta(X)$, and the authors
	raised the question of whether this upper bound could be significantly
	improved:
	
	\begin{problem}
	\label{prob:main}
			Is there a polynomial $P(x)$ such that for every bipartite graph
			$G$ with parts $X$ and $Y$, and
			maximum degree $\Delta$, $\chi'_{int}(G,X) \leq P(\Delta)$?
	\end{problem}

	Asratian (see e.g. \cite{AsratianHaggkvist})
	proved that if a bipartite graph $G$ with parts $X$ and $Y$
	satisfies $d_G(x) \geq d_G(y)$ for all edges $xy \in E(G)$,
	where $x \in X$ and $y \in Y$, and where $d_G(x)$
	denotes the degree of $x$ in $G$, then
	$G$ has an $X$-interval coloring such that each vertex
	$x \in X$ receives colors
	$1,\dots, d_G(x)$ on its incident edges. Thus, for such bipartite graphs,
	the answer to this question is positive.

	In this short note, we prove that Problem \ref{prob:main} has 
	a positive answer for
	every bipartite graph; in fact we show that a cubic polynomial suffices.
	We prove this in Section 2, and in Section 3 we deduce some results 
	on $\chi'_{int}(G,X)$ for particular families of bipartite graphs.


\section{A general upper bound on $\chi'_{int}(G,X)$}

In this paper, ``graphs'' do not have any multiple edges, unless
otherwise stated, while a ``hypergraph'' may have
repeated edges.

We shall denote by $G=(X,Y;E)$ a bipartite graph $G$ with parts $X$ and $Y$,
and edge set $E=E(G)$; $\Delta(X)$ denotes the maximum degree of the vertices in $X$.
 If, in addition, $G$ is $(a,b)$-biregular, then
we shall assume that the vertices of $X$ have degree $a$.

The {\em chromatic index} $\chi'(H)$ of a hypergraph $H$ is the smallest number of colors 
needed for a proper edge coloring of $H$; that is, an edge coloring
where every pair of edges that share a common vertex get different colors.

For a proper edge coloring $\varphi$ of a graph $G$ and $v \in V(G)$ 
we denote by $\varphi(v)$ the set of
colors appearing on edges incident to $v$; $\varphi(v)$ is called
the {\it palette}
of $v$.
 We say that $\varphi$ {\em is interval at $v$}
if $\varphi(v)$ is an interval of integers.

\begin{lemma}
\label{lem:edgehyp}
	If $H$ is a hypergraph with maximum degree $\Delta$ where every edge contains
	at most $k$ vertices,
	then $\chi'(H) \leq k (\Delta-1) +1$.
\end{lemma}
	This follows by applying a greedy edge coloring algorithm, since every edge
	intersects at most $k(\Delta-1)$ other edges.
	
	Denote by $N_G(x)$ the neighborhood of a vertex $x$ in a graph $G$.
	
\begin{lemma}
\label{lem:decomp}
	If $G=(X,Y;E)$ is $(a,b)$-biregular,
	then there is a decomposition of $G$ into at most $b$ subgraphs $F_1,\dots, F_b$,
	where $\Delta(F_i) = a$ and every vertex of $X$ is in precisely
	one of the subgraphs $F_1,\dots, F_b$.
\end{lemma}
\begin{proof}
		From $G$ we form a hypergraph $H$ with vertex set $Y$
		by for every vertex $x$ of $X$ including a hyperedge $N_G(x)$ in $H$. 
		Then $H$ has maximum degree $b$, and every edge in $H$ 
		contains $a$ vertices, so
		by Lemma
		\ref{lem:edgehyp}, there is a proper edge coloring of $H$ with at most 
		$ab$ colors. For each $i = 1,\dots,b$, let $H_i$ be the (edge-induced)
		subhypergraph
		induced by all edges with colors from $\{(i-1)a+1,\dots, (i-1)a +a\}$.
		Then each $H_i$ is a hypergraph with maximum degree at most $a$,
		and so, each $H_i$ corresponds to a subgraph $F_i$ of $G$ with the required properties.
\end{proof}

\begin{proposition}
\label{prop:result}
	If $G =(X,Y;E)$ is $(a,b)$-biregular,
	then it has an $X$-interval edge coloring with at most $ab$ colors.
\end{proposition}

\begin{proof}
	By the preceding lemma, $G$ has a decomposition into at most 
	$b$ subgraphs $F_1,\dots, F_b$ such that $\Delta(F_i) \leq a$ and every vertex
	of $X$ is in precisely one of the subgraphs $F_1,\dots, F_b$. Thus, 
	by K\"onig's edge coloring theorem, each $F_i$ has a proper $a$-edge coloring,
	and this coloring is an $X_{F_i}$-interval coloring of $F_i$, where $X_{F_i} \subseteq X$
	is one part of $F_i$. Now, by using disjoint color
	sets for such colorings of the distinct $F_i$, we obtain the required coloring of
	$G$.
\end{proof}

\begin{theorem}
\label{th:main}
		If $G =(X,Y;E)$ is a bipartite graph, then
		$\chi'_{int}(G,X) \leq \frac{(\Delta(G))^2 (\Delta(G)+1)}{2}$.
\end{theorem}

\begin{proof}
	Let $G$ be a bipartite graph with parts $X$ and $Y$ 
	and $d \in \{1,\dots, \Delta(G)\}$.
	The subgraph $H$
	induced by all edges that are incident with
	vertices of degree $d$ in $X$
	is a subgraph of a $(d,\Delta(G))$-biregular graph, so by 
	Proposition \ref{prop:result},
	$H$ has an $X_H$-interval coloring with at most $d \Delta(G)$ colors,
	where $X_H \subseteq X$ is one part of $H$.
	
	Since every possible vertex degree $d$ in $G$ is in the set 
	$\{1,\dots, \Delta(G)\}$, we deduce that
	$$\chi'_{int}(G,X) \leq
	\sum_{1}^{\Delta(G)} d \Delta(G)=
	\frac{(\Delta(G))^2 (\Delta(G)+1)}{2}.$$
\end{proof}

\begin{remark}
	The  technique in the proof of Lemma \ref{lem:decomp} fails
	if we allow multiple edges in $G$, so Theorem \ref{th:main} is not valid
	for multigraphs. Nevertheless, by proceeding as in the proof of Lemma
	\ref{lem:decomp}, it is straightforward that an
	$(a,b)$-biregular multigraph $G=(X,Y;E)$
	has a decomposition into at most $ab$ subgraphs
	$F_1, \dots, F_{ab}$, where $\Delta(F_i) =a$ and every vertex in $X$
	is in precisely one of these subgraphs. Since K\"onig's edge coloring
	theorem applies to multigraphs, we get a polynomial upper bound on
	$\chi'(G,X)$ for bipartite multigraphs as well, but this time the polynomial
	is of degree $4$.
\end{remark}


\section{Some further bounds}

As first noted in \cite{Avesani}, questions on interval edge colorings
are related to problems on the {\em palette index}
$\check s(G)$ of $G$, defined as the the minimum number
of distinct palettes occurring in a proper edge coloring of $G$.
This notion was introduced quite recently by Hornak et al
\cite{HornakKalinowskiMeszkaWozniak}.

We note the following relation between $\check s(G)$ and $\chi'_{int}(G,X)$.

\begin{proposition}
\label{prop:palette}
	If $G=(X,Y;E)$ is a bipartite graph, then 
	$\chi'_{int}(G,X) \leq \Delta(X) \check s(G)$.
\end{proposition}
\begin{proof}
	Let $\varphi$ be proper edge coloring of $G$ with $\check s(G)$ palettes.
	For any specific palette that appears on a vertex of $X$, we consider the subgraph
	induced by the edges incident with such vertices in $X$. 
	Let $H$ be such a subgraph and let $X_H \subseteq X$ and $Y_H$ be its parts.
	Then all vertices in $X_H$ have the same degree, and all vertices of $Y_H$ have smaller degree in $H$
	than the vertices of $X_H$. Hence, $H$ is
	$\Delta(H)$-edge colorable, and such a coloring of $H$ is an $X_H$-interval
	coloring with $\Delta(H)$ colors.
	Since at most $\check s(G)$ different palettes appear at vertices in $G$, 
	the upper bound follows.
\end{proof}

In \cite{CasselgrenToft3}, it was proved that if a bipartite graph $G=(X,Y;E)$ satisfies
$\Delta(X)=2$, then $\chi'_{int}(G,X) \leq \Delta(G)+1$; moreover,
if $\Delta(X) \geq \Delta(G)-1$, then $\chi'_{int}(G,X) \leq 2 \Delta(G)-2$.

By proceeding as in the proof of Theorem \ref{th:main}, using the aforementioned upper bounds
and combining Proposition \ref{prop:palette} with
Propositions 4.7, 4.9 and Corollary 4.5 in \cite{CasselgrenPetrosyan},
we can deduce that a quadratic upper bound holds for some bipartite graphs.

\begin{corollary}
\label{cor:improve}
	Let $G=(X,Y;E)$ be a bipartite graph. There is a constant
	$C$ such that if
	
	\begin{itemize}
	
	\item[(i)] $\Delta(X) \leq 4$, or
	
	\item[(ii)] $\Delta(G)$ is even and $\Delta(X) \geq \Delta(G)-2$, then 
	
	\end{itemize}

$\chi'_{int}(G,X) \leq C (\Delta(G))^2$.

\end{corollary}
	
	For the special case of bipartite graphs where 
	all vertices in $X$
	have degree $3$,
	a better upper bound 
	was obtained in \cite{Renman}.
	Moreover,
	in view of Corollary \ref{cor:improve}, 
	it is natural to ask the following.
	
	\begin{problem}
	\label{prob:gen}
		Is it true that there is a constant $C$ such that for any bipartite graph $G=(X,Y;E)$,
		 $\chi'_{int}(G,X) \leq C (\Delta(G))^2$?
	\end{problem}
		
	Finally, let us deduce an upper bound of $\chi'_{int}(G,X)$
	for bipartite graphs with small maximum degree.
	In \cite{CasselgrenToft3}, the following upper bounds were obtained:

\begin{itemize}

	\item if $\Delta(G) \leq 4$, then $\chi'_{int}(G,X) \leq 6$;
	
	\item if $\Delta(G) \leq 5$, then $\chi'_{int}(G,X) \leq 15$;

	\item if $\Delta(G) \leq 6$, then $\chi'_{int}(G,X) \leq 33$.
		
\end{itemize}

Here we shall prove the following improvement of the last bound.

\begin{proposition}
\label{prop:smalldegrees}
	If $G=(X,Y;E)$ is a bipartite graph with $\Delta(G) =6$,
	then $\chi'_{int}(G,X) \leq 17$.
\end{proposition}

We shall need the following lemma.

\begin{lemma}
\label{lem:46}
	If $G=(X,Y;E)$ is a bipartite graph with $\Delta(G) =6$ where all vertex degrees
	in $X$ are in $\{1,2,4,5,6\}$, then $\chi'_{int}(G,X) \leq 10$.
\end{lemma}
\begin{proof}
	From $G$ and a copy $G'$ of the graph $G$ we form a $6$-regular 
	multigraph $H$ in the following way:
	
	\begin{itemize}
	
		\item	for each vertex $v$ of odd degree in $G$, we add an edge between
		$v$ and its corresponding vertex $v'$ in $G'$; thereafter we add
		$3-(d_G(v)+1)/2 $ loops at $v$ and $v'$;
	
		\item for each vertex $v$ of even degree in $G$, we add $3-d_G(v)/2$ loops at $v$ and
		at the corresponding vertex $v'$ in $G'$.

	\end{itemize}
	By Petersen's $2$-factor theorem, $H$ can be decomposed into
	three $2$-factors, which in $G$ corresponds to 
	subgraphs $F_1,F_2,F_3$ of maximum degree at most $2$.
	By coloring the edges of $F_i$ alternately by colors $2i-1, 2i$, we obtain
	a proper $6$-edge coloring of $H$. Let $\varphi$ be the restriction of
	this coloring to $G$. Then $\varphi$ is interval at every vertex $v$ of
	$X$ except if 
	$$\varphi(v) \in \{ \{1,2,5,6\}, \{1,2,3,4,6\}, \{1,2,3,5,6\}, \{1,2,4,5,6\},
	\{1,3,4,5,6\}
	\}.$$
	By recoloring some appropriate edges $e$ incident with such vertices by the color 
	$\varphi(e) +6$
	if $\varphi(e) \in \{1,2,3,4\}$, we obtain an $X$-interval coloring of $G$ with $10$ colors.
\end{proof}

\begin{proof}[Proof of Proposition \ref{prop:smalldegrees}]
	By proceeding as in the proof of Theorem \ref{th:main},
	the proposition now follows by combining Lemma \ref{lem:46}
	with the fact that every $(3,6)$-biregular graph has an
	interval $7$-coloring \cite{CasselgrenToft}.
\end{proof}

As a final remark, let us note that Lemma \ref{lem:46} applies to multigraphs, and since this
also holds for the result on interval colorings of $(3,6)$-biregular graphs \cite{CasselgrenToft},
Proposition \ref{prop:smalldegrees} is valid for multigraphs as well.

\section{Acknowledgement}
Carl Johan Casselgren was supported by a grant from the Swedish
Research Council (2017-05077).\\
The author wants to thank Armen S. Asratian for helpful comments and suggestions.

\end{document}